\documentclass[reqno,12pt,a4paper]{amsart}

\usepackage{mathrsfs}
\usepackage{amssymb}
\usepackage{amsfonts}
\usepackage{latexsym}
\usepackage{amsthm}

\usepackage{graphicx}
\def\lb{\label}

\newcommand{\er}[1]{\textrm{(\ref{#1})}}

\begin{document}


\renewcommand{\theequation}{\arabic{section}.\arabic{equation}}
\theoremstyle{plain}
\newtheorem{theorem}{\bf Theorem}[section]
\newtheorem{lemma}[theorem]{\bf Lemma}
\newtheorem{corollary}[theorem]{\bf Corollary}
\newtheorem{proposition}[theorem]{\bf Proposition}
\newtheorem{definition}[theorem]{\bf Definition}
\newtheorem{remark}[theorem]{\it Remark}

\def\a{\alpha}  \def\cA{{\mathcal A}}     \def\bA{{\bf A}}  \def\mA{{\mathscr A}}
\def\b{\beta}   \def\cB{{\mathcal B}}     \def\bB{{\bf B}}  \def\mB{{\mathscr B}}
\def\g{\gamma}  \def\cC{{\mathcal C}}     \def\bC{{\bf C}}  \def\mC{{\mathscr C}}
\def\G{\Gamma}  \def\cD{{\mathcal D}}     \def\bD{{\bf D}}  \def\mD{{\mathscr D}}
\def\d{\delta}  \def\cE{{\mathcal E}}     \def\bE{{\bf E}}  \def\mE{{\mathscr E}}
\def\D{\Delta}  \def\cF{{\mathcal F}}     \def\bF{{\bf F}}  \def\mF{{\mathscr F}}
\def\c{\chi}    \def\cG{{\mathcal G}}     \def\bG{{\bf G}}  \def\mG{{\mathscr G}}
\def\z{\zeta}   \def\cH{{\mathcal H}}     \def\bH{{\bf H}}  \def\mH{{\mathscr H}}
\def\e{\eta}    \def\cI{{\mathcal I}}     \def\bI{{\bf I}}  \def\mI{{\mathscr I}}
\def\p{\psi}    \def\cJ{{\mathcal J}}     \def\bJ{{\bf J}}  \def\mJ{{\mathscr J}}
\def\vT{\Theta} \def\cK{{\mathcal K}}     \def\bK{{\bf K}}  \def\mK{{\mathscr K}}
\def\k{\kappa}  \def\cL{{\mathcal L}}     \def\bL{{\bf L}}  \def\mL{{\mathscr L}}
\def\l{\lambda} \def\cM{{\mathcal M}}     \def\bM{{\bf M}}  \def\mM{{\mathscr M}}
\def\L{\Lambda} \def\cN{{\mathcal N}}     \def\bN{{\bf N}}  \def\mN{{\mathscr N}}
\def\m{\mu}     \def\cO{{\mathcal O}}     \def\bO{{\bf O}}  \def\mO{{\mathscr O}}
\def\n{\nu}     \def\cP{{\mathcal P}}     \def\bP{{\bf P}}  \def\mP{{\mathscr P}}
\def\r{\rho}    \def\cQ{{\mathcal Q}}     \def\bQ{{\bf Q}}  \def\mQ{{\mathscr Q}}
\def\s{\sigma}  \def\cR{{\mathcal R}}     \def\bR{{\bf R}}  \def\mR{{\mathscr R}}
\def\S{\Sigma}  \def\cS{{\mathcal S}}     \def\bS{{\bf S}}  \def\mS{{\mathscr S}}
\def\t{\tau}    \def\cT{{\mathcal T}}     \def\bT{{\bf T}}  \def\mT{{\mathscr T}}
\def\f{\phi}    \def\cU{{\mathcal U}}     \def\bU{{\bf U}}  \def\mU{{\mathscr U}}
\def\F{\Phi}    \def\cV{{\mathcal V}}     \def\bV{{\bf V}}  \def\mV{{\mathscr V}}
\def\P{\Psi}    \def\cW{{\mathcal W}}     \def\bW{{\bf W}}  \def\mW{{\mathscr W}}
\def\o{\omega}  \def\cX{{\mathcal X}}     \def\bX{{\bf X}}  \def\mX{{\mathscr X}}
\def\x{\xi}     \def\cY{{\mathcal Y}}     \def\bY{{\bf Y}}  \def\mY{{\mathscr Y}}
\def\X{\Xi}     \def\cZ{{\mathcal Z}}     \def\bZ{{\bf Z}}  \def\mZ{{\mathscr Z}}
\def\O{\Omega}

\newcommand{\gA}{\mathfrak{A}}
\newcommand{\gB}{\mathfrak{B}}
\newcommand{\gC}{\mathfrak{C}}
\newcommand{\gD}{\mathfrak{D}}
\newcommand{\gE}{\mathfrak{E}}
\newcommand{\gF}{\mathfrak{F}}
\newcommand{\gG}{\mathfrak{G}}
\newcommand{\gH}{\mathfrak{H}}
\newcommand{\gI}{\mathfrak{I}}
\newcommand{\gJ}{\mathfrak{J}}
\newcommand{\gK}{\mathfrak{K}}
\newcommand{\gL}{\mathfrak{L}}
\newcommand{\gM}{\mathfrak{M}}
\newcommand{\gN}{\mathfrak{N}}
\newcommand{\gO}{\mathfrak{O}}
\newcommand{\gP}{\mathfrak{P}}
\newcommand{\gQ}{\mathfrak{Q}}
\newcommand{\gR}{\mathfrak{R}}
\newcommand{\gS}{\mathfrak{S}}
\newcommand{\gT}{\mathfrak{T}}
\newcommand{\gU}{\mathfrak{U}}
\newcommand{\gV}{\mathfrak{V}}
\newcommand{\gW}{\mathfrak{W}}
\newcommand{\gX}{\mathfrak{X}}
\newcommand{\gY}{\mathfrak{Y}}
\newcommand{\gZ}{\mathfrak{Z}}

\def\ve{\varepsilon}   \def\vt{\vartheta}    \def\vp{\varphi}    \def\vk{\varkappa}

\def\Z{{\mathbb Z}}    \def\R{{\mathbb R}}   \def\C{{\mathbb C}}    \def\K{{\mathbb K}}
\def\T{{\mathbb T}}    \def\N{{\mathbb N}}   \def\dD{{\mathbb D}}


\def\la{\leftarrow}              \def\ra{\rightarrow}            \def\Ra{\Rightarrow}
\def\ua{\uparrow}                \def\da{\downarrow}
\def\lra{\leftrightarrow}        \def\Lra{\Leftrightarrow}


\def\lt{\biggl}                  \def\rt{\biggr}
\def\ol{\overline}               \def\wt{\widetilde}
\def\no{\noindent}


\let\ge\geqslant                 \let\le\leqslant
\def\lan{\langle}                \def\ran{\rangle}
\def\/{\over}                    \def\iy{\infty}
\def\sm{\setminus}               \def\es{\emptyset}
\def\ss{\subset}                 \def\ts{\times}
\def\pa{\partial}                \def\os{\oplus}
\def\om{\ominus}                 \def\ev{\equiv}
\def\iint{\int\!\!\!\int}        \def\iintt{\mathop{\int\!\!\int\!\!\dots\!\!\int}\limits}
\def\el2{\ell^{\,2}}             \def\1{1\!\!1}
\def\sh{\sharp}
\def\wh{\widehat}
\def\bs{\backslash}

\def\all{\mathop{\mathrm{all}}\nolimits}
\def\Area{\mathop{\mathrm{Area}}\nolimits}
\def\arg{\mathop{\mathrm{arg}}\nolimits}
\def\const{\mathop{\mathrm{const}}\nolimits}
\def\det{\mathop{\mathrm{det}}\nolimits}
\def\diag{\mathop{\mathrm{diag}}\nolimits}
\def\diam{\mathop{\mathrm{diam}}\nolimits}
\def\dim{\mathop{\mathrm{dim}}\nolimits}
\def\dist{\mathop{\mathrm{dist}}\nolimits}
\def\Im{\mathop{\mathrm{Im}}\nolimits}
\def\Iso{\mathop{\mathrm{Iso}}\nolimits}
\def\Ker{\mathop{\mathrm{Ker}}\nolimits}
\def\Lip{\mathop{\mathrm{Lip}}\nolimits}
\def\rank{\mathop{\mathrm{rank}}\limits}
\def\Ran{\mathop{\mathrm{Ran}}\nolimits}
\def\Re{\mathop{\mathrm{Re}}\nolimits}
\def\Res{\mathop{\mathrm{Res}}\nolimits}
\def\res{\mathop{\mathrm{res}}\limits}
\def\sign{\mathop{\mathrm{sign}}\nolimits}
\def\span{\mathop{\mathrm{span}}\nolimits}
\def\supp{\mathop{\mathrm{supp}}\nolimits}
\def\Tr{\mathop{\mathrm{Tr}}\nolimits}
\def\BBox{\hspace{1mm}\vrule height6pt width5.5pt depth0pt \hspace{6pt}}
\def\where{\mathop{\mathrm{where}}\nolimits}
\def\as{\mathop{\mathrm{as}}\nolimits}


\newcommand\nh[2]{\widehat{#1}\vphantom{#1}^{(#2)}}
\def\dia{\diamond}

\def\Oplus{\bigoplus\nolimits}



\def\qqq{\qquad}
\def\qq{\quad}
\let\ge\geqslant
\let\le\leqslant
\let\geq\geqslant
\let\leq\leqslant
\newcommand{\ca}{\begin{cases}}
\newcommand{\ac}{\end{cases}}
\newcommand{\ma}{\begin{pmatrix}}
\newcommand{\am}{\end{pmatrix}}
\renewcommand{\[}{\begin{equation}}
\renewcommand{\]}{\end{equation}}
\def\eq{\begin{equation}}
\def\qe{\end{equation}}
\def\[{\begin{equation}}
\def\bu{\bullet}

\title[{Estimates of length of spectrum}]
        {Sharp spectral estimates for periodic matrix-valued Jacobi operators}
\date{\today}

\def\Wr{\mathop{\rm Wr}\nolimits}
\def\BBox{\hspace{1mm}\vrule height6pt width5.5pt depth0pt \hspace{6pt}}

\def\Diag{\mathop{\rm Diag}\nolimits}

\date{\today}
\author[Anton Kutsenko]{Anton Kutsenko}
\address{Laboratoire de M\'ecanique Physique, UMR CNRS 5469,
Universit\'e Bordeaux 1, Talence 33405, France,  \qqq email \
kucenkoa@rambler.ru }

\subjclass{81Q10 (34L40 47E05 47N50)} \keywords{matrix-valued Jacobi
operator, Jacobi matrix, spectral estimates, measure of spectrum}

\maketitle

\begin{abstract}
For the periodic matrix-valued Jacobi operator $J$ we obtain the
estimate of the Lebesgue measure of the spectrum $|\s(J)|\le4
\min_n\Tr(a_na_n^*)^\frac12$, where $a_n$ are off-diagonal elements
of $J$. Moreover estimates of width of spectral bands are obtained.
\end{abstract}

\section{Introduction}
\setcounter{equation}{0}

We consider a self-adjoint matrix-valued Jacobi operator $
J:\ell^2(\Z)^m\to\ell^2(\Z)^m$ given by
\[\lb{001}
 (Jy)_n=a_n y_{n+1}+b_ny_n+a_{n-1}^*y_{n-1},\ \ n\in\Z,\ \
 y_n\in\C^m,\ \ y=(y_n)_{n\in\Z}\in\ell^2(\Z)^m,
\]
where $a_n$ and $b_n=b_n^*$ are $p$-periodic sequences of the
complex $m\ts m$ matrices. It is well known (see e.g. \cite{CG},
\cite{CGR}, \cite{KKu}, \cite{KKu1}) that the spectrum of this
operator $\s(J)=\s_{ac}(J)\cup\s_{p}(J)$, where absolutely
continuous part $\s_{ac}(J)$ is a union of finite number of
intervals and $\s_{p}(J)$ consists of finite number of eigenvalues
of infinite multiplicity. Note that if $\det a_n\ne0$ for all
$n=1,...,p$, then $\s(J)=\s_{ac}(J)$ and $\s_{p}(J)=\es$ always. The
main goal of this paper is to obtain estimate of length of spectrum
$\s(J)$. We don't know such estimates for the matrix-valued Jacobi
operators.

\begin{theorem}\lb{T1}
The Lebesgue measure of spectrum of $J$ satisfy the following
estimate
\[\lb{002}
 {\rm mes}(\s(J))\le4
\min_n\Tr(a_na_n^*)^\frac12.
\]
\end{theorem}
{\no\bf Remark 1.}Note that this estimate does not depend on period
$p$ and coefficients $b_n$, i.e. changing $p$ and $b_n$ we can't
sufficiently increase the length of spectrum of $J$. If $a_n=0$ for
some $n\in\Z$, then this estimate is sharp.

{\no\bf 2. ($m=1$)} For the scalar case $m=1$ we have the estimate
(see e.g. \cite{DS}, \cite{Ku}, \cite{KKr})
\[\lb{002a}
 {\rm mes}(\s(J))\le4|a_1a_2...a_p|^{1\/p},
\]
We reach equality for the case of discrete Shr\"odinger operator
$J^0$ with $a_n^0=1$, $b_n^0=0$. Estimate \er{002} is better than
\er{002a}, since $\min|a_n|\le|a_1a_2...a_p|^{1\/p}$.

{\no\bf 3.} The result similar to \er{002a} was obtained in
\cite{PR} for general non periodic scalar case ($m=1$).

{\no \bf 4.} From the Proof of Theorem \ref{T1} we obtain estimates
of spectral bands, see \er{008}.

{\no\bf Example  (sharpness).} We construct the Jacobi matrix $J$
whose spectrum satisfy ${\rm mes}(\s(J))=4\min_n(\Tr
a_na_n^*)^{\frac12}$. Let $J$ be Jacobi matrix with elements
$a_n=I_m$ ($m\ts m$ identical matrix) and $b_n=\diag(4k)_{k=1}^m$
for any $n$. Since all $a_n$ and $b_n$ are diagonal matrix, then $J$
is unitarily equivalent to the direct sum of scalar Jacobi
operators. In our case this is the direct sum of shifted discrete
Shr\"odinger operators $\os_{k=1}^m(J^{0}+4kI)$ ($I$ is identical
operator). Then
$$
 \s(
 J)=\bigcup_{k=1}^{m}\s(J^{0}+4kI)=\bigcup_{k=1}^{m}[-2+4k,2+4k]=[2,2+4m],
$$
which gives us ${\rm mes}(\s(J))=4m=4(\Tr a_na_n^*)^{\frac12}$.

\section{Proof of Theorem \ref{T1}}

Without lost of generality we may assume that
$\min_n\Tr(a_na_n^*)^{\frac12}=\Tr(a_0a_0^*)^{\frac12}$ and $p\ge3$.
It is well known (see e.g. \cite{KKu}, \cite{KKu1}) that $J$ is
unitarily equivalent to the operator
$\cJ=\int_{[0,2\pi)}^{\os}K(x)dx$ acting in
$\int_{[0,2\pi)}^{\os}\cH dx$, where $\cH=\C^{pm}$ and $pm\ts pm$
matrix $K(x)$ is given by
\[\lb{004}
 K(x)=\left(\begin{array}{ccccc} b_1 & a_1 & 0 & ... & e^{-ix}a_0^* \\
                                      a_1^* & b_2 & a_2 & ... & 0 \\
                                      0 & a_2^* & b_3 & ... & 0 \\
                                      ... & ... & ... & ... & ... \\
                                      e^{ix} a_0 & 0 & 0 & ... & b_p
  \end{array}\right)=K_0+K_1(x),
\]
where
\[\lb{005}
 K_0=\left(\begin{array}{ccccc} b_1 & a_1 & 0 & ... & 0 \\
                                      a_1^* & b_2 & a_2 & ... & 0 \\
                                      0 & a_2^* & b_3 & ... & 0 \\
                                      ... & ... & ... & ... & ... \\
                                      0 & 0 & 0 & ... & b_p
  \end{array}\right),\ \ K_1(x)=\left(\begin{array}{ccccc} 0 & 0 & 0 & ... & e^{-ix}a_0^* \\
                                      0 & 0 & 0 & ... & 0 \\
                                      0 & 0 & 0 & ... & 0 \\
                                      ... & ... & ... & ... & ... \\
                                      e^{ix} a_0 & 0 & 0 & ... & 0
  \end{array}\right).
\]
The spectrum $\s(J)$ is
\[\lb{005a}
 \s(J)=\bigcup_{x\in[0,2\pi]}\s(K(x)).
\]
From \er{004}-\er{005} we obtain
\[\lb{006}
 K_0-|K_1|\le K(x)\le K_0+|K_1|,\ \ x\in[0,2\pi]
\]
where
\[\lb{007}
 |K_1|=(K_1K_1^*)^{\frac12}=\left(\begin{array}{ccccc} (a_0^*a_0)^{\frac12} & 0 & 0 & ... & 0 \\
                                      0 & 0 & 0 & ... & 0 \\
                                      0 & 0 & 0 & ... & 0 \\
                                      ... & ... & ... & ... & ... \\
                                      0 & 0 & 0 & ... &
                                      (a_0a_0^*)^{\frac12}
  \end{array}\right)
\]
does not depend on $x$. Let $\l_1(x)\le...\le\l_N(x)$ be eigenvalues
of $K(x)$ and let $\l_1^{\pm}\le...\le\l_N^{\pm}$ be eigenvalues of
$K_0\pm|K_1|$. Using \er{006} we obtain
\[\lb{008}
 \l^-_n\le\l_n(x)\le\l^+_n,\ \ x\in[0,2\pi]
\]
which with \er{005a} gives us
\[\lb{009}
 \s(J)=\bigcup_{x\in[0,2\pi]}\{\l_n(x)\}_{n=1}^N\ss\bigcup_{n=1}^N[\l_n^-,\l_n^+].
\]
Then
\[\lb{010}
 {\rm
 mes}(\s(J))\le\sum_{n=1}^N(\l_n^+-\l_n^-)=2\Tr|K_1|=4\Tr(a_na_n^*)
 \BBox
\]

{\it Acknowledgements.} I would to express thanks to prof. E.
Korotyaev for useful discussions and remarks. Also I want to thank
prof. B. Simon for useful comments and refferences to the paper
\cite{PR}.

\no Many thanks to Michael J. Gruber, who told me that instead of
$\|a_n\|\rank a_n$ (which was in the first version of this paper) is
better to use $\Tr(a_na_n^*)^{\frac12}$ in \er{002}.


\begin{thebibliography} {9999}
\setlength{\itemsep}{-\parskip} \footnotesize

\bibitem [CG]{CG}  Clark, S.; Gesztesy, F. On Weyl–-Titchmarsh theory for
singular finite difference Hamiltonian systems, J. Comput. Appl.
Math., 171 (2004) 151–184.

\bibitem [CGR]{CGR} Clark, S.; Gesztesy, F.; Renger, W.
Trace formulas and Borg-type theorems for matrix-valued Jacobi and
Dirac finite difference operators. J. Diff. Eq. 219 (2005),
144--182.

\bibitem[DS] {DS} P. Deift, B. Simon. Almost periodic Schr\"odinger operators III.
The absolutely continuous spectrum in one dimension. Commun. Math.
Phys., 90, 389–411 (1983).

\bibitem[K] {K} Kato T. Perturbation Theory for Linear Operators. Springer (February 15,
1995).

\bibitem [Ku] {Ku} Kutsenko A. Estimates of Parameters for Conformal Mappings Related to a Periodic Jacobi
Matrix. Journal of Mathematical Sciences,     Volume 134, Number 4 /
April 2006, Pages 2295-2304.

\bibitem [KKr]{KKr}  Korotyaev, E.; Krasovsky, I.
Spectral estimates for periodic Jacobi matrices, Commun. Math. Phys.
234(2003), 517-532.


\bibitem [KKu]{KKu} Korotyaev, E., Kutsenko, A.
    Lyapunov functions for periodic matrix-valued Jacobi operators,
    AMS translations Series 2, 225 (2008), 117—-131.

\bibitem[KKu1] {KKu1} Korotyaev, E., Kutsenko, A. Borg type uniqueness Theorems for
periodic Jacobi operators with matrix valued coefficients. Proc. of
the AMS, Volume 137, Number 6, June 2009, Pages 1989–-1996.

\bibitem[L] {L} Y. Last. On the measure of gaps and spectra for discrete 1D
Schr\"odinger operators. Commun. Math. Phys., 149, 347–-360 (1992).

\bibitem[PR]{PR} A. Poltoratski, C. Remling. Reflectionless Herglotz Functions and Jacobi
Matrices, Commun. Math. Phys. Volume 288 Number 3(2009), 1007--1021.

\bibitem[RS]{RS} M. Reed ; B. Simon. Methods of modern mathematical physics. IV.
Analysis of operators. Academic Press, New York-London, 1978.

\bibitem[S]{S} B. Simon, Orthogonal polynomials on the unit circle, Part 1 and Part 2, AMS, Providence, RI, 2005.

\bibitem[S1]{S1} B. Simon, Trace Ideals and Their Applications: Second
Edition. Mathematical Surveys and Monographs vol. 120, 2005.

\bibitem[Te] {Te}  G. Teschl, Jacobi Operators and Completely Integrable Nonlinear
Lattices, Mathematical Surveys and Monographs, vol. 72, AMS, Rhode
Island, 2000.

\bibitem [vM]{vM}  P. van Moerbeke. The spectrum of Jacobi matrices.
Invent. Math. 37 (1976), no. 1, 45--81.

\end{thebibliography}
\end{document}